\documentclass[11pt,a4paper,english,reqno]{amsart}
\usepackage{amsmath,amssymb,amsfonts,epsfig,mathrsfs}
\usepackage[T1]{fontenc}

\usepackage{color}
\usepackage{array}
\usepackage{amsthm}
\usepackage{amstext}
\usepackage{graphicx}
\usepackage{setspace}
\usepackage[margin=2.5cm]{geometry}
\usepackage{bbm}
\usepackage{color}
\usepackage{enumitem}
\allowdisplaybreaks[4]

\usepackage{amscd,psfrag}
\usepackage{yhmath}
\usepackage[mathscr]{eucal}

\usepackage{slashed}

\makeatletter
\pdfpageheight\paperheight
\pdfpagewidth\paperwidth

\usepackage{epstopdf}
\usepackage{chngcntr}
\counterwithin{figure}{section}
\usepackage{mathrsfs}

\usepackage{indentfirst}	

\usepackage[normalem]{ulem}
\theoremstyle{plain}

\newtheorem*{theorem*}{Theorem}

\newtheorem*{remark*}{Remark}
\newtheorem*{sideremark*}{Side Remark}\newtheorem*{mt*}{Main Theorem}

\newtheorem*{claim*}{Claim}
\newtheorem*{q*}{Question}

\newtheorem*{lemma*}{Lemma}

\newtheorem*{corollary*}{Corollary}
\newtheorem*{proposition*}{Proposition}

\newcommand{\R}{\mathbb{R}}

\newcommand{\na}{\nabla}
\newcommand{\dd}{{\rm d}}
\newcommand{\p}{\partial}
\newcommand{\e}{\epsilon}

\newcommand{\map}{\rightarrow}
\newcommand{\G}{\Gamma}

\newcommand{\E}{{\mathscr{E}}}

\allowdisplaybreaks[4]

\def\XXint#1#2#3{{\setbox0=\hbox{$#1{#2#3}{\int}$ }
\vcenter{\hbox{$#2#3$ }}\kern-.6\wd0}}

\newcommand{\mres}{\mathbin{\vrule height 1.6ex depth 0pt width
0.13ex\vrule height 0.13ex depth 0pt width 1.3ex}}

\newcommand{\str}{{\mathbb{S}}}

\numberwithin{equation}{section}
\numberwithin{figure}{section}

\title{Regularity of Desingularised Models for Vortex Filaments in Incompressible Viscous Flows: \\A Geometrical Approach}

\author{Siran Li}
\address{Siran Li: Department of Mathematics, Rice University, MS 136
P.O. Box 1892, Houston, Texas, 77251, USA.}

\email{\texttt{Siran.Li@rice.edu}}



\date{\today}

\begin{document}

\maketitle

\begin{center}
\emph{Dedicated to Bob Hardt, with Admiration and Friendship}
\end{center}

\begin{abstract}
We establish the regularity of weak solutions for the vorticity equation associated to a family of desingularised models for vortex filament dynamics in 3D incompressible viscous flows. These include and generalise the classical model ``of an allowance for the thickness of the vortices'' due to Louis Rosenhead in 1930. Our approach is based on an interplay between the geometry of vorticity and analytic inequalities in Sobolev spaces.
\end{abstract}

\section{Introduction}
\subsection{The PDE for vorticity}
In this paper we study the regularity issues of vortex filament dynamics in an incompressible viscous fluid in $\R^3$. The  \emph{vorticity} $\omega$ is a vectorfield of  physical significance: it measures the ``size of rotations'' in the fluid and  plays a central r\^{o}le in the regularity theory of the fluid flow ({\it cf.} Wolibner \cite{w} and Yudovich \cite{y}). 

The dynamics of $\omega$ is described by the following PDE:
\begin{equation}\label{vorticity eq}
\frac{\p\omega}{\p t} + (u\cdot \na ) \omega - \nu \Delta \omega = \str\cdot \omega\qquad \text{ in } [0,T^\star[\times\R^3.
\end{equation}
Here, $\omega: [0,T^\star[\times\R^3\map\R^3$ is the vorticity, $u: [0,T^\star[\times\R^3\map\R^3$ is the velocity, the constant $\nu>0$ is the kinematic viscosity, and $\str=\str(u):[0,T^\star[\times\R^3\map\R^3\otimes\R^3$ is the rate-of-strain tensor of the fluid, defined as
\begin{equation}\label{str}
\str := \frac{1}{2}\big(\na u + \na^\top u\big).
\end{equation} 
Eq.~\eqref{vorticity eq} is supplemented by the incompressibility condition
\begin{equation}\label{incompressible}
{\rm div}\, u = 0\qquad \text{ in } [0,T^\star[\times\R^3
\end{equation}
and the initial condition
\begin{equation}\label{id}
\omega|\{t=0\}=\omega_0\qquad \text{ on } \{0\} \times \R^3.
\end{equation}

\subsection{Vortex filaments}\label{subsec: filaments} Our work is concerned with models of line vortices and vortex filaments. Related physical notions and mathematical models are summarised in \S \S \ref{subsec: filaments}--\ref{subsec: rosenhead}, which are essentially taken from \S 2 in Berselli--Gubinelli \cite{bg}. 

A \emph{line vortex} is a singular distribution in which infinite vorticity is concentrated on a curve $\gamma$, such that the circulation $\G>0$ around a closed circuit threaded by $\gamma$ is finite. $\G$ is called the \emph{strength} of the line vortex. One may view a line vortex as obtained via the limiting process of pinching a \emph{vortex filament} --- that is, a vortex tube surrounded by the fluid --- to the curve $\gamma$, with the strength being kept constant ({\it cf.} Helmholtz \cite{h}).

We consider the case that the curve supporting the vorticity is a \emph{knot}, {\it i.e.}, a smooth simple closed curve $\gamma: [0,1] \map \R^3$ with $\gamma(0)=\gamma(1)$. Throughout $\gamma$ shall be identified with its image.

Formally, the vorticity vectorfield $\omega$ is given by Eq.~(2.4) in \cite{bg}:
\begin{equation}\label{formal eq}
\omega(t,x) = \G \int_0^1 \delta\big(x-\gamma(t,\xi)\big) \gamma_\xi(t,\xi)\,\dd \xi\qquad \text{ for each } x \in \R^3,\, t \in [0,T^\star[.
\end{equation}
Here $\delta$ is the Dirac delta function, $\xi \in [0,1]$ is the arclength parameter, and $\gamma_\xi = \p\gamma\slash\p\xi$. In measure-theoretic notations we write $\omega(t,x)=\G\int_0^1 \gamma_\xi(t,\xi)\,\dd\mu(\xi)$, with the measure $\mu$ being the restriction of the $1$-dimensional Hausdorff measure to the curve: $\mu \equiv \mathcal{H}^1 \mres \{x=\gamma(t,\bullet)\}$.

\subsection{Biot--Savart laws}\label{subsec: biot--savart}

As is well-known, for the Navier--Stokes and Euler equations, the vorticity $\omega$ is related to the velocity $u$ of the fluid via
\begin{equation}
\omega = {\rm curl}\, u = \na \times u.
\end{equation}
In $\R^3$, one may represent the operator ${\rm curl}^{-1}$ as a singular integral:
\begin{equation}\label{bs, simple}
u(t,x) = -\frac{1}{4\pi} \int_{\R^3} \frac{x-y}{|x-y|^3}\times \omega(t,y)\,\dd y.
\end{equation}
Such a relation is known as a \emph{Biot--Savart law}.

Consider the case of a line vortex. Assuming that the knot $\gamma$ is transported by the velocity vectorfield, one may deduce from Eqs.\,\eqref{formal eq} and \eqref{bs, simple} that
\begin{equation}\label{vortex eq}
\frac{\p\gamma}{\p t}(t,\xi) = -\frac{\G}{4\pi} \int_0^1 \frac{\gamma(t,\xi)-\gamma(t,\eta)}{|\gamma(t,\xi)-\gamma(t,\eta)|^3} \times \gamma_\eta(t,\eta)\,\dd\eta;
\end{equation}
see \cite{bg, s, bb1}. Near the diagonal $\{\xi=\eta\} \subset [0,1]\times[0,1]$ the PDE~\eqref{vortex eq} is highly singular.

\subsection{Rosenhead approximation}\label{subsec: rosenhead}

In 1930, Rosenhead \cite{r} proposed and analytically studied a desingularised model for Eq.~\eqref{vortex eq}. The paper \cite{r} begins as such:
\begin{quotation}
This paper is an attempt to investigate the effect on the configuration of vortices in the wake behind a cylinder of an allowance for the thickness of the vortices...
\end{quotation}

The idea of Rosenhead's approximation is to smear out the singularity of Eq.~\eqref{vortex eq} on the diagonal by considering the desingularised vortex equation:
\begin{equation}\label{regularised vortex eq}
\frac{\p\gamma}{\p t}(t,\xi) = -\frac{\G}{4\pi} \int_0^1 \frac{\gamma(t,\xi)-\gamma(t,\eta)}{\Big[\big(\gamma(t,\xi)-\gamma(t,\eta)\big)^2+\mu^2\Big]^{3\slash 2}} \times \gamma_\eta(t,\eta)\,\dd\eta
\end{equation}
for some constant $\mu>0$. See also Moore \cite{m} for an application of this model in  numerical computations for aircraft trailing vortices.

In effect, one may view Rosenhead's desingularised model~\eqref{regularised vortex eq} as obtained via a \emph{modified Biot--Savart law}. As in \cite{bg, bb1}, Eq.\,\eqref{regularised vortex eq} amounts to expressing $u$ in terms of $\omega$ by
\begin{align}\label{bs-rosenhead}
u(t,x) = -\frac{1}{4\pi}\int_{\R^3} \na\phi(x-y)\times\omega(t,y)\,\dd y,
\end{align}
with the potential $\phi:[0,T^\star[\times\R^3\map\R$ given by
\begin{equation}\label{potential-rosenhead}
\phi(z) = \frac{\G}{\sqrt{|z|^2+\mu^2}}.
\end{equation}
Note that $\phi(z)$ becomes completely regular as $|z| \map 0$. 

One should compare Rosenhead's model with the usual Biot--Savart law for $u={\rm curl}^{-1}\omega$, namely Eq.\,\eqref{bs, simple}. The latter equation can be obtained  from Eq.~\eqref{potential-rosenhead} by setting $\mu=0$.

\subsection{Partial desingularisations}\label{subsec: Partial desingularisations} Our main goal is to investigate ``partially  desingularised'' models for Eq.\,\eqref{vortex eq}. We shall analyse the Biot--Savart law~\eqref{bs-rosenhead} with the potential
\begin{equation}\label{potential--general}
\phi_\delta(z) = \frac{\G}{\sqrt{|z|^2+\mu^2|z|^\delta}}\qquad \text{ for a positive parameter $\delta$}.
\end{equation}
In the case $0<\delta<2$, we obtain a desingularised model which is \emph{more singular} than Rosenhead's approximation ($\delta=0$) in Eqs.~\eqref{bs-rosenhead}\eqref{potential-rosenhead}\eqref{regularised vortex eq}, but \emph{less singular} than $u={\rm curl}^{-1}\omega$ ($\delta=2$).

Using methods pioneered by Constantin--Fefferman \cite{cf} in the study of geometrical regularity criteria for the Navier--Stokes equations, we  establish the regularity of weak solutions for the vorticity equation~\eqref{vorticity eq}, under the partially desingularised Biot--Savart laws~\eqref{potential--general} with
\begin{equation*}
0\leq \delta <1.
\end{equation*}

The precise statement of our results are given in \S \ref{sec: result}.

\subsection{Related works}
For the background on the mathematical analysis of fluid dynamical PDEs, we refer to Constantin--Foias \cite{cf}, Temam \cite{temam}, Ladyzhenskaya \cite{l}, Lemari\'{e}-Rieusset \cite{lr}, Seregin \cite{s}, Galdi \cite{g}, Robinson--Rodrigo--Sadowski \cite{rrs}, and many others.

The dynamics of vortices is an important topic in aero- and hydro-dynamics; see Saffman \cite{s} and Chorin \cite{c} for a comprehensive treatment. The existence, uniqueness, and stability properties of various PDE models for vortex dynamics have been studied; {\it cf.} Banica--Vega \cite{bv}, Aiki--Iguchi \cite{ai}, Jerrard--Seis \cite{js}, Lions--Majda \cite{lm}, etc. The study of vortex dynamics in incompressible \emph{viscous} fluid flows attracts much attention in recent works; \textit{cf.} Enciso--Luc\`{a}--Peralta-Salas \cite{elp}. 

Regarding the Rosenhead model, a rigorous analytic study was first carried out by Berselli--Bessaih \cite{bb1}. 
It is extended to more general stochastic contexts by Bessaih--Gubinelli--Russo \cite{bgr}; see also Flandoli \cite{f}.

The global existence of the smooth solution for Eqs.~\eqref{vorticity eq}\eqref{potential--general} was proved by Berselli--Gubinelli \cite{bg} under a few mild assumptions on the \emph{Fourier side} of the potential $\phi$ in the Biot--Savart law; see ``Hypothesis A'' on p.698 therein. The results in \cite{bg} cover a wide range of desingularised models for vortex filament dynamics with modified Biot--Savart laws, including the Rosenhead approximation. In this paper, we are concerned with desingularised models of a different nature: the desingularisation effects take place on the \emph{physical side} of $\phi$. The potentials considered in \S \ref{subsec: Partial desingularisations} (see Eq.~\eqref{potential--general}) appear to be still more singular at  the singularity $z=0$.

The key new feature of this paper is to apply ideas and techniques from the works on  \emph{geometric regularity criteria} for the Navier--Stokes equations to analyse the desingularised models for vortex filament dynamics. First proposed by Constantin--Fefferman \cite{cf}, the geometric regularity criteria can be summarised as follows (here the vorticity and velocity are related by the usual Biot--Savart law \eqref{bs, simple}): ``If for all $t \in [0,T^\star[$ the angle between the vorticity vectorfield of a weak solution for the Navier--Stokes equations at nearby points satisfies certain uniform H\"{o}lder conditions in space, then it is automatically strong up to the time $T^\star$''.  Such criteria have been further developed in many works; see, {\it e.g.},  \cite{bb2, bb3, dg, g, gr, gm, z, ty1, ty2, v, ct, ber, bei,  me1, me2, bp, zwz}. This list is by no means exhaustive.

Finally, intriguing linkages between topological--differential geometrical works on \emph{knot energies} and analytic studies on line vortex dynamics deserve further explorations. Such linkages are suggested by Eq.\,\eqref{vortex eq}. We refer to Freedman--He--Wang \cite{fhw} and O'Hara \cite{o} for knot energies.

\section{Main Result}\label{sec: result}
\subsection{Theorem}\label{thm}The main result of our paper is as follows. For the convenience of readers, we reproduce Eqs.~\eqref{vorticity eq}\eqref{bs-rosenhead}\eqref{potential--general} in the statement below.

\begin{quotation}
Let $\omega \in L^\infty(0,T^\star; L^1 \cap W^{-1,2}(\R^3;\R^3))\cap L^2(0,T^\star; L^2(\R^3;\R^3))$ be a weak solution for the vorticity equation:
\begin{equation*}
\frac{\p \omega}{\p t} + (u\cdot\na) \omega -\nu\Delta \omega  = \str\cdot\omega
\qquad \text{ in } [0,T^\star[\times\R^3.
\end{equation*}
Assume that $u$ is related to $\omega$ via the modified Biot--Savart law:
\begin{equation*}
u(t,x) = -\frac{1}{4\pi}\int_{\R^3} \na\phi_\delta(x-y)\times\omega(t,y)\,\dd y
\end{equation*}
where, for positive constants $\G$ and $\mu$, 
\begin{equation*}
\phi_\delta(z) = \frac{\G}{\sqrt{|z|^2+\mu^2|z|^\delta}}.
\end{equation*}
Furthermore, assume that $\{\omega(t,\bullet)\}_{t \in [0,T^\star[}$ is supported on a compact set in $\R^3$. Then, the vorticity $\omega$ has higher regularity as follows: $\omega \in L^\infty(0,T^\star; L^2(\R^3;\R^3)) \cap L^2(0,T^\star;W^{1,2}(\R^3;\R^3))$ whenever $\delta \in [0,1[$.
\end{quotation}

\subsection{Remark}
The assumption $\omega \in L^\infty(0,T^\star; L^1(\R^3;\R^3))$ is due to the uniform boundedness of total circulation. We put
\begin{equation}\label{Sigma}
\Sigma_0 := \sup_{0\leq t <T^\star} \int_{\R^3}|\omega(t,x)|\,\dd x < \infty.
\end{equation} 
The compact support for $\{\omega(t,\bullet)\}_{t \in [0,T^\star[}$ means that the vortex filament does not become infinitely large up to time $T^\star$. These agree with the discussions in \S \ref{subsec: filaments} on the physical model. 

In this paper, weak solutions are understood in the distributional sense as usual; {\it e.g.}, for the vorticity equation~\eqref{vorticity eq}, one needs to integrate against arbitrary test functions of the form $\psi(t)\chi(x)$ for $\phi \in C^\infty_0(]-1,T^\star[)$ and $\chi \in C^\infty_c(\R^3)$. 

\subsection{Strategy of the proof}
We represent the rate-of-strain  tensor $\str=\str(u)$ as a singular integral of $\omega$ via the modified Biot--Savart law. Then, we show that the vorticity stretching term
\begin{equation}\label{stretch}
\mathfrak{S}(t) := \int_{\R^3} \str(t,x) : \omega(t,x) \otimes \omega(t,x)\,\dd x
\end{equation}
can be controlled by the \emph{enstrophy}
\begin{equation}\label{enstrophy}
\E(t) := \frac{1}{2}\int_{\R^3} |\omega(t,x)|^2\,\dd x.
\end{equation}
(Throughout, for $3\times 3$ matrices $\mathfrak{n}$ and $\mathfrak{p}$ we write $\mathfrak{n}:\mathfrak{p}\equiv {\rm tr}(\mathfrak{n}\cdot\mathfrak{p})$.)
This is done by exploring the r\^{o}le of the angle $\angle(\omega(t,x), \omega(t,y))$ played in the singular integral as in \cite{cf}. 


\section{Inequalities}
In this section we summarise a few well-known analytic inequalities for Sobolev functions that shall be utilised in \S \ref{sec: proof}. 

\subsection{Interpolation for $L^p$}\label{subsec: interpolation, lp}

Let $p_0$, $p_1$ be such that $0<p_0<p_1\leq \infty$. For any $0 \leq \theta \leq 1$  define $p_\theta$ by $\frac{1}{p_\theta} = \frac{1-\theta}{p_0} + \frac{\theta}{p_1}$. Then, for any $n=1,2,3,\ldots$ and any $f \in L^{p_0}(\R^n) \cap L^{p_1}(\R^n)$, there holds
\begin{equation*}
\|f\|_{L^{p_\theta}(\R^n)} \leq \|f\|_{L^{p_0}(\R^n)}^{1-\theta} \|f\|^\theta_{L^{p_1}(\R^n)}.
\end{equation*}

\subsection{Gagliardo--Nirenberg--Sobolev inequality} \label{subsec: gns}
Let $f:\R^n\map \R$. Fix $1\leq q,r<\infty$ and $m=1,2,3,\ldots$. Suppose that $\alpha \in \R$ and $j \in \mathbb{N}$ satisfy 
\begin{align*}
\frac{j}{m} \leq \alpha \leq 1,\qquad \frac{1}{p} = \frac{j}{n} + \Big(\frac{1}{r}-\frac{m}{n}\Big)\alpha + \frac{1-\alpha}{q}. 
\end{align*}
Then there exists a constant $C$ depending only on $m$, $n$, $j$, $q$, $r$, and $\alpha$ such that 
\begin{equation*}
\|D^j f\|_{L^p(\R^n)} \leq C \|D^m f\|^\alpha_{L^r(\R^n)} \|f\|^{1-\alpha}_{L^q(\R^n)}.
\end{equation*}

\subsection{Hardy--Littlewood--Sobolev inequality}\label{subsec: hls}
Let $f\in L^p(\R^n)$ and $g \in L^s(\R^n)$ with $1<p,s<\infty$. Assume that $0<\lambda<n$ satisfies $1/p + 1/s+\lambda/n=2$. Then there is a constant $C$ depending only on $p$, $n$, and $\lambda$ such that 
\begin{equation*}
\Bigg|\iint_{\R^n\times\R^n} f(x)|x-y|^{-\lambda} g(y)  \,\dd x\,\dd y\Bigg| \leq C \|f\|_{L^p(\R^n)} \|g\|_{L^s(\R^n)}.
\end{equation*}

\section{Proof of \ref{thm} Theorem}\label{sec: proof}
\subsection{Preliminary energy estimate}\label{subsec: energy est}
Multiplying $\omega$ to both sides of the vorticity equation~\eqref{vorticity eq} and integrating over space, we obtain
\begin{equation*}
\frac{\dd\E}{\dd t}(t) + \nu \int_{\R^3} |\na\omega(t,x)|^2\,\dd x = \mathfrak{S}(t).
\end{equation*}
The right-hand side is the vorticity stretching term given by Eq.~\eqref{stretch}. One should note that this evolution equation for $\E$ is understood in the sense of distributions over $[0,T^\star[$.

\subsection{Singular integral representation of $\str$}\label{subsec: lemma} Now we show
\begin{lemma*}
In the setting of Theorem~\ref{thm}, the rate-of-strain tensor $\str$ can be represented as
\begin{align}
\str(x) &= {\rm p.v.}\,\int_{\R^3} \bigg\{ \frac{\G\delta(2-\delta)}{4} \mu^2 A^{-\frac{3}{2}}(|x-y|) \cdot |x-y|^{\delta-4} + \frac{3\G}{8} B^2(|x-y|) A^{-\frac{5}{2}}(|x-y|) \bigg\} \cdot\nonumber\\
&\qquad\qquad \cdot \bigg\{ (x-y) \times \omega(y) \otimes (x-y) + (x-y) \otimes (x-y) \times \omega(y) \bigg\} \,\dd y.
\end{align}
Here $A, B: [0,\infty[\map[0,\infty[$ are given by
\begin{eqnarray}
&&A(r) := r^2 + \mu^2 r^\delta,\label{A}\\
&& B(r) := 2+\delta\mu^2 r^{\delta-2}.\label{B}
\end{eqnarray}
\end{lemma*}

Most estimates in this paper hold pointwise in $t$; so, unless otherwise declared, we always suppress the $t$-variable. Throughout $p.v.$ denotes the principal value of singular integrals; for vectors $a=(a^1,a^2,a^3)$ and $b=(b^1,b^2,b^3)\in\R^3$, $a \otimes b$ is the $3\times 3$ matrix $\{a^ib^j\}_{1\leq i,j\leq 3}$. 

\subsection{Proof of \ref{subsec: lemma} Lemma} First let us show that
\begin{equation}\label{xx}
\na u(x) = {\rm p.v.} \int_{\R^3} \na\na\phi_\delta(x-y) \times \omega(y)\,\dd y,
\end{equation}
where
\begin{equation*}
(\na\na\phi_\delta\times\omega)^j_i \equiv \na_i (\na\phi_\delta\times \omega)^j
\end{equation*}
for $i,j\in\{1,2,3\}$.  This is an equality for $3\times 3$ matrices.

Indeed, direct computation yields that
\begin{align}\label{yy}
\na \phi_\delta(z) &= -\frac{\G}{2} \na\Big(|z|^2+\mu^2|z|^\delta\Big)\Big(|z|^2+\mu^2|z|^\delta\Big)^{-\frac{3}{2}}\nonumber\\
&=-\frac{\G}{2} \Big(2z + \delta \mu^2|z|^{\delta-2}z\Big)\Big(|z|^2+\mu^2|z|^\delta\Big)^{-\frac{3}{2}}.
\end{align}
It is locally integrable on $\R^3$. This enables us to compute the weak Hessian $\na\na\phi_\delta$ via the dominated convergence theorem as follows.

Take an arbitrary test function $\chi \in C^\infty_c(\R^3)$. The above paragraph ensures that
\begin{align*}
-\langle\chi, \na \na \phi_\delta\rangle = \lim_{\e \searrow 0} \int_{\{|x|\geq\e\}} \na\phi_\delta(x)\otimes\na\chi(x)\,\dd x,
\end{align*}
where the left-hand side is the  paring of a distribution with a  test function. In view of integration by parts and the divergence theorem, it further equals
\begin{equation*}
\lim_{\e \searrow 0} \Bigg\{ -\int_{\{|x|\geq\e\}} \na\na\phi_\delta(x)\chi(x)\,\dd x + \int_{\{|x|=\e\}}\chi(x) \na\phi_\delta(x)\otimes\frac{x}{|x|}\,\dd \mathcal{H}^2(x) \Bigg\}.
\end{equation*}
For the second term, a change of variable leads to
\begin{align*}
\int_{\{|x|=\e\}}\chi(x) \na\phi_\delta(x)\otimes\frac{x}{|x|}\,\dd \mathcal{H}^2(x) = \e^2 \int_{\{|x|=1\}} \chi(\e x) \na\phi_\delta(\e x) \otimes x \,\dd\mathcal{H}^2(x).
\end{align*}
By the definition of $\phi_\delta$ in Eq.\,\eqref{potential--general}, the right-hand side is controlled by $C(\mu,\G)\e^{2-\frac{3\delta}{2}}\|\chi\|_{L^\infty(\R^3)}$. For any $\delta \in [0,1[$ this tends to zero as $\e\searrow 0$. Thus Eq.~\eqref{xx} follows.

The previous arguments justify the computation of $\na\na\phi_\delta$ by directly taking $\na$ to the final line in Eq.~\eqref{yy}. In local coordinates, we get
\begin{align*}
\na_i\na_j\phi_\delta(z) &= -\frac{\G}{2}\na_i\Bigg\{ \frac{2z^j + \delta\mu^2|z|^{\delta-2}z^j}{A^{\frac{3}{2}}(|z|)} \Bigg\}\\
&= -\frac{\G}{2} A^{-3}(|z|) \Bigg\{A^{\frac{3}{2}}(|z|) \bigg[ 2 \tilde{\delta}_{ij} +\delta \mu^2 \tilde{\delta}_{ij}|z|^{\delta-2} + \delta(\delta-2)\mu^2|z|^{\delta-4}z^iz^j \bigg] \\
&\qquad\qquad\qquad-\frac{3}{2}\Big(2z^j + \delta\mu^2z^j|z|^{\delta-2} \Big(A^{\frac{1}{2}}(|z|)\na_i A(|z|)\Big) \Bigg\}.
\end{align*}
Here $\tilde{\delta}_{ij}$ is the Kronecker delta symbol  (lest one confuses it with the parameter $\delta \in [0,1[$). 

Next, note the simple identity
\begin{equation*}
\na_i A(|z|) = z^i B(|z|),
\end{equation*}
from which we infer that
\begin{align*}
\na\na\phi_\delta(z) &= -\frac{\G}{2} A^{-\frac{3}{2}}(|z|) A^{-\frac{3}{2}}(|z|) \Big[\delta\mu^2\tilde{\delta}+\delta(\delta-2)\mu^2|z|^{\delta-4}z\otimes z\Big] \\
&\qquad+ \frac{3\G}{4}A^{-\frac{5}{2}}(|z|) B^2(|z|) z\otimes z.
\end{align*}
Thus, for some anti-symmetric matrix $\mathfrak{m} \in \mathfrak{so}(3,\R)$, there holds
\begin{align}\label{zz}
&\na\na\phi_\delta(x-y)\times\omega(y) \nonumber\\
&\qquad= \mathfrak{m} + \bigg\{(x-y)\times\omega(y)\otimes(x-y)\bigg\}\cdot\nonumber\\
&\qquad\qquad\cdot \bigg\{ -\frac{\G}{2}\delta(\delta-2)\mu^2A^{-\frac{3}{2}}(|x-y|) |x-y|^{\delta-4} + \frac{3\G}{4}A^{-\frac{5}{2}}(|x-y|) B^2(|x-y|) \bigg\}.
\end{align}

We substitute Eq.~\eqref{zz} into Eq.~\eqref{xx} to get the  singular integral representation for $\na u$. The lemma in \S \ref{subsec: lemma} follows immediately by  symmetrising the resulting expression.

\subsection{Vorticity stretching term $\mathfrak{S}$}\label{subsec: lemma'}
The singular integral representation for $\str$ in \S \ref{subsec: lemma} implies
\begin{lemma*}
The vorticity stretching term $\mathfrak{S}=\int_{\R^3}\str:\omega\otimes\omega\,\dd x$ can be bounded as follows:
\begin{equation}\label{S-lemma}
|\mathfrak{S}| \leq \int_{\R^3} |\omega(x)|^2 \Bigg\{{\rm p.v.}\int_{\R^3}\Big[K^{(1)}(|x-y|)+K^{(2)}(|x-y|)\Big]|\omega(y)|\,\dd y\Bigg\}\,\dd x,
\end{equation}
where 
\begin{eqnarray}
&& K^{(1)}(r) := \frac{\G \delta(2-\delta)}{4} \mu^2 r^{\delta-2} A^{-\frac{3}{2}}(|z|), \label{K1} \\
&& K^{(2)}(r) := \frac{3\G}{8} r^2 B^2(r) A^{-\frac{5}{2}}(r). \label{K2}
\end{eqnarray}
\end{lemma*}

Recall that $A$ and $B$ are given in \S \ref{subsec: lemma}. Thanks to this lemma, the control for $\mathfrak{S}$ reduces to estimating the double integral of the product of $|\omega(x)|^2$, $|\omega(y)|$, and a kernel $\lesssim \mathcal{O}(|x-y|^\lambda)$ for some exponent $\lambda>0$. This is tackled by  the Hardy--Littlewood--Sobolev inequality.

\subsection{Proof of \ref{subsec: lemma'} Lemma} In \S \ref{subsec: lemma} we proved that
\begin{align*}
\str &= {\rm p.v.} \int_{\R^3} \Big[K^{(1)}(|x-y|)+K^{(2)}(|x-y|)\Big] \cdot\nonumber\\
&\qquad\qquad \cdot \bigg\{ (x-y) \times \omega(y) \otimes (x-y) + (x-y) \otimes (x-y) \times \omega(y) \bigg\} \,\dd y.
\end{align*}
Following the crucial observations in Constantin--Fefferman \cite{cf}, by writing
\begin{equation*}
\widehat{z} := \frac{z}{|z|} \qquad \text{ for any } z \in \R^3,
\end{equation*}
one may express the vorticity stretching term as
\begin{align*}
\mathfrak{S} &= \int_{\R^3} |\omega(x)|^2 \Bigg\{{\rm p.v.}\int_{\R^3} |\omega(y)|\Big[K^{(1)}(|x-y|)+K^{(2)}(|x-y|)\Big]\cdot\\
&\qquad\qquad\qquad\qquad\cdot\begin{bmatrix}
\widehat{x-y} \times \widehat{\omega}(y) \otimes \widehat{x-y}\\
+\\
\widehat{x-y}\otimes \widehat{x-y}\times \widehat{\omega}(y)
\end{bmatrix} : \Big[\widehat{\omega}(x) \otimes  \widehat{\omega}(x)\Big] \Bigg\}\,\dd x.
\end{align*}
But there holds
\begin{align*}
\begin{bmatrix}
\widehat{x-y} \times \widehat{\omega}(y) \otimes \widehat{x-y}\\
+\\
\widehat{x-y}\otimes \widehat{x-y}\times \widehat{\omega}(y)
\end{bmatrix} : \Big[\widehat{\omega}(x) \otimes  \widehat{\omega}(x)\Big] = \mathscr{D}\Big(\widehat{x-y}, \widehat{\omega}(x), \widehat{\omega}(y)\Big),
\end{align*}
where for arbitrary unit column vectors $e_1,e_2,e_3\in\R^3$ we write
\begin{equation*}
\mathscr{D}(e_1,e_2,e_3) := e_1\cdot e_3 \det(e_1|e_2|e_3).
\end{equation*}
Thus elementary calculus gives us the pointwise estimate:
\begin{align}\label{calc}
\Big|\mathscr{D}\Big(\widehat{x-y}, \widehat{\omega}(x), \widehat{\omega}(y)\Big) \Big| \leq \Big|\sin\angle\Big(\omega(t,x),\omega(t,y)\Big)\Big|.
\end{align}
We can now conclude \ref{subsec: lemma'} by 
bounding this term na\"{i}vely by $1$. 


\subsection{Conclusion of the proof}

We resume from estimate~\eqref{S-lemma} for the vorticity stretching term $\mathfrak{S}$. Thanks to Eqs.~\eqref{A}\eqref{B}\eqref{K1} and \eqref{K2}, pointwise we have
\begin{eqnarray}
&&K^{(1)}(|z|) \leq \frac{\G \delta (2-\delta)}{4\mu} |z|^{-2-\frac{\delta}{2}},\label{K est 1}\\
&& K^{(2)}(|z|) \leq 
\frac{3\G}{4\mu^{5}}|z|^{-2-\frac{\delta}{2}} \qquad \text{ if } |z| \leq \bigg( \frac{\delta\mu^2}{2} \bigg)^{\frac{1}{2-\delta}},\label{K est 2}\\
&& K^{(2)}(|z|) \leq \frac{3\G}{\mu^5} |z|^{2-\frac{5\delta}{2}}\qquad \text{ if } |z| \geq \bigg( \frac{\delta\mu^2}{2} \bigg)^{\frac{1}{2-\delta}}.\label{K est 3}
\end{eqnarray}

Let us also introduce the constants
\begin{eqnarray}
&& \eta := \max\bigg\{ \frac{3\G}{\mu^5}, \, \frac{\G \delta (2-\delta)}{4\mu} \bigg\},\label{const1}\\
&& \Lambda := {\rm diam}^{2-\frac{5\delta}{2}}  \bigg({\rm spt}\,\{\omega(t,\bullet)\}_{t \in [0,T^\star[}\bigg)\label{const2}.
\end{eqnarray}
Applying bounds~\eqref{K est 1}\eqref{K est 2} and \eqref{K est 3} on the kernels to \eqref{S-lemma} in \S \ref{subsec: lemma'}, we deduce that
\begin{align}\label{aaa}
|\mathfrak{S}| \leq \eta\Lambda \|\omega\|^2_{L^2} \|\omega\|_{L^1} + \eta \iint_{\R^3 \times \R^3} |\omega(x)|^2 |\omega(y)| |x-y|^{-2-\frac{\delta
}{2}} \,\dd x\,\dd y.
\end{align} 

Now we invoke the Hardy--Littlewood--Sobolev inequality in \ref{subsec: gns}. We take $p=1+\e$, $n=3$, and $s=2+\frac{\delta}{2}$ therein, for some $\e>0$ to be determined. Then 
\begin{align}\label{hls, application}
\Bigg|\iint_{\R^3 \times \R^3} |\omega(x)|^2 |\omega(y)| |x-y|^{-2-\frac{\delta
}{2}} \,\dd x\,\dd y\Bigg| \leq C_1 \|\omega\|^2_{L^{2(1+\e)}(\R^3)} \|\omega\|_{L^s(\R^3)},
\end{align}
where $C_1$ depends only on $\delta$ and $\e$, with the index $s = \frac{6(1+\e)}{2(1+4\e)-\delta(1+\e)}$.

Recall that $\delta \in [0,1[$; hence, by additionally requiring that $\e \in ]0,\frac{1}{2}[$, we have $s \in ] \frac{6}{4-\delta}, \frac{6}{2-\delta}[$. So we can fix an $\e$ (depending only on $\delta$) in $]0,\frac{1}{2}[$ once and for all to warrant that $1<s<2$. Consequently, one may consider $s$ being fixed from now on. As a result, the constant $C_1$ in the inequality~\eqref{hls, application} depends only on $\delta$. Then the interpolation in \ref{subsec: interpolation, lp} gives us
\begin{align}\label{lp interpolation-application}
\|\omega\|_{L^s(\R^3)} \leq \|\omega\|_{L^1(\R^3)}^{\frac{2}{s}-1} \|\omega\|_{L^2(\R^3)}^{2-\frac{2}{s}}.
\end{align}

On the other hand, the Gagliardo--Nirenberg--Sobolev inequality in \ref{subsec: gns} implies that 
\begin{align}\label{gns-application}
\|\omega\|^2_{L^{2(1+\e)}(\R^3)} \leq C_2 \|\na \omega\|_{L^2(\R^3)}^{\frac{3\e}{1+\e}} \|\omega\|_{L^2(\R^3)}^{\frac{2-\e}{1+\e}}.
\end{align}
Here $C_2$ depends only on $\e$, hence only on $\delta$ as in the previous paragraph.

Putting together Eqs.~\eqref{hls, application}\eqref{lp interpolation-application} and \eqref{gns-application}, one obtains 
\begin{align*}
\Bigg|\iint_{\R^3 \times \R^3} |\omega(x)|^2 |\omega(y)| |x-y|^{-2-\frac{\delta
}{2}} \,\dd x\,\dd y\Bigg| \leq C_3 \|\na \omega\|_{L^2(\R^3)}^{\frac{3\e}{1+\e}} \|\omega\|_{L^2(\R^3)}^{2-\frac{2}{s}+\frac{2-\e}{1+\e}},
\end{align*}
where $C_3=C_1C_2 (\Sigma_0)^{\frac{2}{s}-1}$ depends only on $\delta$ and the finite total circulation $\Sigma_0$ (see Eq.~\eqref{Sigma}).

In addition, the simple  inequality $ab \leq \frac{a^p}{p} + \frac{b^q}{q}$ for $\frac{1}{p}+\frac{1}{q}=1$, $1<p<\infty$, and $a,b\geq 0$ implies that
\begin{align}\label{young, app}
\Bigg|\iint_{\R^3 \times \R^3} |\omega(x)|^2 |\omega(y)| |x-y|^{-2-\frac{\delta
}{2}} \,\dd x\,\dd y\Bigg| \leq \frac{\nu}{2} \|\na \omega\|^2_{L^2(\R^3)} + C_4 \|\omega\|_{L^2(\R^3)}^{2\kappa},
\end{align}
where the constant
\begin{equation}\label{const3}
C_4 = \frac{2-\e}{2(1+\e)} \frac{(C_3)^{\frac{2(1+\e)}{3\e}}}{\big(\frac{\nu}{2}\big)^{\frac{3\e}{2-\e}}},
\end{equation}
and the exponent
\begin{equation*}
\kappa = \frac{(1+\e)}{2-\e} \Big(2-\frac{2}{s} + \frac{2-\e}{1+\e}\Big).
\end{equation*}
Since $s \in ]1,2[$ and $\e \in ]0,\frac{1}{2}[$ (thus $\frac{1+\e}{2-\e} \in ]\frac{1}{2}, 1[$), it is crucial to see that
\begin{equation}\label{kappa}
\frac{1}{4}<\kappa<2,
\end{equation}
where $\kappa$ depends only on $\delta$.

To conclude, we substitute Eqs.~\eqref{young, app}\eqref{aaa} into the energy estimate in \S \ref{subsec: energy est} to get
\begin{align}\label{diff ineq}
\frac{\dd \E}{\dd t} + \frac{\nu}{2} \|\na \omega\|^2_{L^2(\R^3)} \leq  \eta\Lambda \Sigma_0 \E(t) + \eta C_4 \E^\kappa
\end{align}
for the enstrophy $\E(t)=\int_{\R^3}|\omega(t,x)|^2\,\dd x$. But the hypotheses of the theorem require that $\E \in L^2([0,T^\star[)$. Hence, in light of Eq.~\eqref{kappa}, the differential inequality~\eqref{diff ineq} implies that $\E \in L^\infty([0,T^\star[)$ and that $\na \omega \in L^2([0,T^\star[ \times \R^3;\R^3 \otimes \R^3)$. This concludes the proof.

\section{Remarks}
A closer examination on the dependence of constants ({\it e.g.}, from Eqs.~\eqref{const1}\eqref{const2}\eqref{const3} and \eqref{diff ineq}) shows that we have obtained an upper bound for 
${\rm ess\,sup}_{t \in [0,T^\star[} \|\omega(t,\bullet)\|_{L^2(\R^3)}^2$ that is proportional to the strength $\G$ of vortex filaments, a positive power of the diameter of the support of $\omega$, the total circulation $\Sigma_0$, a positive power of $T^\star$, and a positive power of the initial enstrophy $\int_{\R^3}|\omega_0|^2\,\dd x$. Also, this upper bound is \emph{inverse proportional} to $\mu^5$ and $\nu^\varsigma$, where $\varsigma = \frac{3\e}{2-\e} \in ]0,1[$, and where the  regularisation parameter is assumed to satisfy $0<\mu\ll 1$. 

Nevertheless, we cannot directly pass to the inviscid limit by sending $\e \searrow 0$ (hence $\varsigma \searrow 0$): it corresponds to the endpoint case of the Hardy--Littlewood--Sobolev inequality ($p=1$ in \ref{subsec: hls}).

Meanwhile, we have obtained an upper bound for $\|\na \omega(t,\bullet)\|^2_{L^2(\R^3)}$, which have the same dependence on all the parameters as for  
${\rm ess\,sup}_{t \in [0,T^\star[} \|\omega(t,\bullet)\|_{L^2(\R^3)}^2$, except that it is inverse proportional to $\nu^{1+\varsigma}$.

\end{document}